\documentclass[amscd,amssymb,11pt,reqno]{amsart}
\def\R{{\mathbb R}}

\newtheorem{Thm}{Theorem}[section]

\newtheorem{Cor}[Thm]{Corollary}
\newtheorem{Lem}[Thm]{Lemma}
\newtheorem{Prop}[Thm]{Proposition}
\newtheorem{Def}[Thm]{Definition}
\newtheorem{Rmk}[Thm]{Remark}

\usepackage{graphicx}
\usepackage{latexsym}
\begin{document}

%\begin{Large}

\vspace{1.5 cm}

%\begin{center}{ Last modified Dec 2, 2004  } \end{center}

\title[The geometry of $L_0$.]
      {The geometry of $L_0$.}

\author{N.J.Kalton, A.Koldobsky, V.Yaskin and M.Yaskina}

\begin{abstract}
Suppose that we have the unit Euclidean ball in $\R^n$ and construct new bodies using
three operations - linear transformations, closure in the radial metric and multiplicative
summation defined by $\|x\|_{K+_0L} = \sqrt{\|x\|_K\|x\|_L}.$ We prove that in dimension 3
this procedure gives all origin symmetric convex bodies, while this is no longer true in
dimensions 4 and higher. We introduce the concept of embedding of a normed space in $L_0$
that naturally extends the corresponding properties of $L_p$-spaces with $p\ne0$, and show
that the procedure described above gives exactly the unit balls of subspaces of $L_0$ in
every dimension. We provide Fourier analytic and geometric characterizations of spaces
embedding in $L_0$, and prove several facts confirming the place of $L_0$ in the scale
of $L_p$-spaces.
\end{abstract}

\maketitle

\section{Introduction}

Suppose that we have the unit Euclidean ball in $\R^n$ and are allowed to
construct new bodies using three operations - linear tranformations,
multiplicative summation and closure in the radial metric. The {\it multiplicative
sum} $K+_0L$ of star bodies $K$ and $L$ is defined by
\begin{equation} \label{0sum}
\|x\|_{K+_0L} = \sqrt{\|x\|_K\|x\|_L},
\end{equation}
where $\|x\|_K=\min\{a\ge 0:\ x\in aK\}$ is the Minkowski functional of a star
body $K.$ What class of bodies do we get from the unit ball by means of
these three operations?

We are going to prove that in dimension $n=3$ we get all origin-symmetric convex bodies, while in
dimension 4 and higher this is no longer the case. However, the class of bodies that we get in
arbitrary dimension also has a clear interpretation. We introduce the concept of embedding in $L_0$
and show that the bodies that we get by means of these three operations are exactly the unit balls
of spaces that embed in $L_0.$

The idea of this interpretation comes from a similar result for $L_p$-spaces with
$p\in [-1,1],\ p\neq 0.$
Namely, if we replace the multiplicative summation by $p$-summation
\begin{equation} \label{psum}
\|x\|_{K+_pL}=\left( \|x\|_{K}^p+\|x\|_L^p \right)^{1/p}
\end{equation}
then we get the unit balls of all spaces that embed in $L_p.$ The case $p=1$ is well-known (see
\cite[Corollary 4.1.12]{G2}) and the unit balls of subspaces of $L_1$ have a clear geometric
meaning - these are the polar projection bodies (see \cite{B}). On the other hand, it was proved by
Goodey and Weil [GW] that if $p=-1$ (this case corresponds to the radial summation) then we get the
class of intersection bodies in $\R^n.$ As shown in \cite{K4}, intersection bodies are the unit
balls of spaces that embed in $L_{-1}.$ The concept of embedding in $L_p,\ p<0$ was introduced in
\cite{K3} as an analytic extension of the same property for $p>0,$ see  \cite{KK} for related
results. The result of Goodey and Weil can easily be extended to $p\in (-1,1), \ p\neq 0.$ Note
that this construction provides a continuous (except for $p=0$) path from polar projection bodies
to intersection bodies, which is important for understanding the duality between projections and
sections of convex bodies. One of the goals of this article is to fill the gap in this scheme at
$p=0$ and better understand the geometry of this intermediate case.

Another interesting similarity of our result with other values of $p$ is that for $p=1$ the
procedure defined above gives all origin-symmetric convex bodies only in dimension 2. This follows
from a result of Schneider \cite{S} that every origin-symmetric convex body is a polar projection
body only in dimension 2. When $p=-1$ we get all origin-symmetric convex bodies only in dimensions
4 and lower, because, by results from \cite{G1}, \cite{Z}, \cite{GKS}, only in these dimensions
every origin-symmetric convex body is an intersection body. The transition between the dimensions 2
and 3 in the case $p=1$ and the transition between the dimensions 4 and 5 in the case $p=-1$
directly correspond to the transition between the affirmative and negative answers in the Shephard
and Busemann-Petty problems, respectively. It would be interesting to find a similar geometric
result corresponding to the transition between dimensions 3 and 4 in the case $p=0$. We refer the
reader to the book \cite[Chapter 6]{K5} for more details and history of the connection between
convex geometry and the theory of $L_p$-spaces.

\bigbreak
\section{The definition of embedding in $L_0$.}\label{sec2}

A compact set $K$ in $\R^n$ is called an origin-symmetric star body if every
straight line passing through the origin crosses the boundary of $K$ at exactly
two points, the boundary is continuous, and the origin is an interior point of $K.$
We denote by $(\R^n,\|\cdot\|_K)$ the Euclidean space equipped with the Minkowski
functional of the body $K.$ Clearly, $(\R^n,\|\cdot\|_K)$ is a normed space if and only
if the body $K$ is convex. Throughout the paper, we write $(\R^n,\|\cdot\|)$ meaning
that $\|\cdot\|$ is the Minkowski functional of some origin-symmetric star body.

A well-known result of P.L\'{e}vy, see \cite[p. 189]{BL} or \cite[Section 6.1]{K5}, is that a space
$(\mathbb{R}^n, \|\cdot\|)$ embeds into $L_p$, $p>0$ if and only if there exists a finite Borel
measure $\mu$ on the unit sphere so that, for every $x\in\mathbb{R}^n$,
\begin{eqnarray}{\label{Def:L_p>0}}
\|x\|^p=\int_{S^{n-1}} |(x, \xi ) |^p d\mu(\xi).
\end{eqnarray}
On the other hand, the definition of embedding in $L_p$ with $p<0$ from \cite{K3} implies that a
space $(\mathbb{R}^n, \|\cdot\|)$ embeds into $L_{p}$, $p \in (-n,0)$ if and only if there exists a
finite symmetric measure $\mu$ on the sphere $S^{n-1}$ so that for every test function $\phi$,
\begin{equation} \label{Def:L_p<0}
\int_{\mathbb{R}^n}\|x\|^{p}\phi(x)dx=\int_{S^{n-1}}d\mu(\xi)
\int_{\mathbb{R}}|t|^{-p-1}\hat{\phi}(t\xi)dt.
\end{equation}
Both representations (\ref{Def:L_p>0}) and (\ref{Def:L_p<0}) are invariant with
respect to $p$-summation.
This gives an idea of defining embedding in $L_0$ by means of a representation that is
invariant with respect to multiplicative summation. Note that the multiplicative summation
is the limiting case of $p$-summation as $p\to 0.$

\begin{Def}\label{Def:L_0}
We say that a space $(\mathbb{R}^n, \|\cdot\|)$ embeds in $L_0$ if there exist a finite Borel
measure $\mu$ on the sphere $S^{n-1}$ and a constant $C\in \R$ so that, for every
$x\in\mathbb{R}^n$,
\begin{equation} \label{logrepr}
\ln \|x\| =\int_{S^{n-1}} \ln |(x, \xi ) | d\mu(\xi) + C.
\end{equation}
\end{Def}
While being similar to (\ref{Def:L_p>0}) and (\ref{Def:L_p<0}), this definition has
its unique features. First, the measure $\mu$ must be a probability measure on $S^{n-1}.$
In fact, put $x=k y$, $k>0$ in (\ref{logrepr}). Then
$$\ln k+ \ln \|y\|=\int_{S^{n-1}} \ln  k\  d\mu(\xi)+\int_{S^{n-1}} \ln |( y, \xi ) |
d\mu(\xi)+C$$
and, again by (\ref{logrepr}) with $x=y$, we get
$\ln k= \int_{S^{n-1}} \ln  k\  d\mu(\xi),$ so
$\int_{S^{n-1}}  d\mu(\xi) =1.$

Secondly, the constant $C$ depends on the norm
and can be computed precisely. In order to compute this constant, integrate
the equality (\ref{logrepr}) over the uniform measure on the unit sphere. We get

\begin{eqnarray*}
C\cdot|S^{n-1}| &=&\int_{S^{n-1}} \ln\|x\| dx - \int_{S^{n-1}}\int_{S^{n-1}}\ln\left|(
x,\theta) \right|d{\mu}(\theta) dx\\
                &=&\int_{S^{n-1}} \ln\|x\| dx - \int_{S^{n-1}}\int_{S^{n-1}}\ln\left|(
x,\theta) \right|dx \, d{\mu}(\theta) \\
&=&\int_{S^{n-1}} \ln\|x\| dx - \int_{S^{n-1}}\ln\left|( x,\theta) \right|dx ,
\end{eqnarray*}
since $\int_{S^{n-1}}\ln\left|( x,\theta) \right|dx$ is rotationally invariant and, therefore, is a
constant for $\theta\in S^{n-1}$, and  $\mu$ is a probability measure.

To compute the latter integral, use the well-known formula (see  \cite[Section 6.4]{K5})
$$\int_{S^{n-1}}\left|( x,\theta)\right|^p dx =
 \frac{2\pi^{(n-1)/2}\Gamma((p+1)/2)}{\Gamma((n+p)/2)}.$$

Differentiating with respect to $p$ and letting $p=0$ we get
$$\int_{S^{n-1}}\ln\left|( x,\theta) \right| dx=
\pi^{(n-1)/2}\left[
\frac{\Gamma'(1/2)}{\Gamma(n/2)}-\sqrt{\pi}\frac{\Gamma'(n/2)}{\Gamma^2(n/2)}\right]$$

Note that
$$|S^{n-1}|= \frac{2\pi^{n/2}}{\Gamma(n/2)},$$
so
$$C= \frac{1}{|S^{n-1}|}\int_{S^{n-1}} \ln\|x\| dx -\frac{1}{2\sqrt{\pi}}\Gamma'(1/2)+
\frac{1}{2}\frac{\Gamma'(n/2)}{\Gamma(n/2)}. $$

Let us remark that Definition \ref{Def:L_0} is equivalent to the
following.  A finite-dimensional normed space $X=(\mathbb
R^n,\|\cdot\|)$ embeds into $L_0$ if and only if there is a
probability space $(\Omega,\mu)$ and a linear map $T:X\to\mathcal
M(\Omega,\mu)$ (where $\mathcal M(\Omega,\mu)$ denotes the space
of $\mu-$measurable functions on $\Omega$) such that
$$ \ln\|x\| =\int_{\Omega}\ln |Tx(\omega)|\,d\mu(\omega),\qquad x\in X.$$
(Here the integrals are also assumed to converge.) Indeed if such
an operator $T$ exists we can write it in the form
$$ Tx(\omega)=h(\omega)(x,\xi(\omega)), \qquad x\in X$$ where
$h:\Omega\to \mathbb R^+$ and $\xi:\Omega\to S^{n-1}$ are
measurable.  Then
$$ \int_{S^{n-1}} \ln|(x,\xi(\omega))|dx>-\infty$$ so that it
follows for some $x\in S^{n-1}$, $\omega\to \ln|(x,\xi(\omega)|$
is $\mu-$integrable.  Hence so is $\ln h$ and further
$$ \ln\|x\| =\int \ln h(\omega) \,d\mu(\omega) +
\int \ln|(x,\xi(\omega))|d\mu(\omega).$$  Now we can induce a
probability measure $\mu'$ on $S^{n-1}$ by
$\mu'(B)=\mu\{\omega:\xi(\omega)\in B\}$ and we have the same
situation as Definition 2.1.

One advantage of this viewpoint is that we can  make sense of the
statement that an infinite-dimensional Banach space embeds into
$L_0$.

\bigbreak
\section{A Fourier analytic characterization of subspaces of $L_0$}

As usual, we denote by ${\mathcal S}(\R^n)$ the space of infinitely differentiable rapidly
decreasing functions on $\R^n$ (test functions), and by ${\mathcal S}^{'}(\R^n)$ the space
of distributions over ${\mathcal S}(\R^n).$

We say that a distribution is positive (negative) outside of the origin in $\R^n$ if it assumes
non-negative (non-positive) values on non-negative Schwartz's test functions with compact support
outside of the origin.

The Fourier transform of a distribution $f$ is defined by
$\langle\hat{f}, {\phi}\rangle= \langle f, \hat{\phi} \rangle$
for every test function $\phi.$

Let $\phi$ be an integrable function on $\mathbb{R}^n$ that is also integrable on hyperplanes, let
$\xi \in S^{n-1}$, and let $t \in \mathbb{R}^n$. Then
$$\mathcal{R}\phi(\xi;t)=\int_{(x,\xi)=t}\phi(x) dx $$
is the {\it Radon transform of $\phi$ in the direction $\xi$ at the point $t$}. A simple connection
between the Fourier and Radon transforms is that for every fixed $\xi \in \mathbb{R}^n \setminus
\{0\}$
\begin{equation} \label{radon-fourier}
\hat{\phi}(s\xi)=\left(\mathcal{R}\phi(\xi;t)\right)^\wedge(s), \ \ \ \forall s \in \mathbb{R}
\end{equation}
where in the right hand side we have the Fourier transform of the function $
t\to\mathcal{R}\phi(\xi;t)$.

The fact that the Fourier transform is useful in the study of subspaces of $L_p$ has been known for
a long time. A well-known result of P.Levy is that a finite dimensional normed space
$(\R^n,\|\cdot\|)$ embeds isometrically in $L_p,\ 0<p\le 2$ if and only if $\exp(-\|\cdot\|^p)$ is
a positive definite function on $\R^n.$ It was proved in \cite{K2} that a space $(\R^n,\|\cdot\|)$
embeds isometrically in $L_p,\ p>0,$ $p \notin 2\mathbb{N}$ if and only if the Fourier transform of
the function $\Gamma(-p/2)\|x\|^{p}$ (in the sense of distributions) is a positive distribution
outside of the origin. If $-n<p<0$ a similar fact was proved in \cite{K3}: a space
$(\R^n,\|\cdot\|)$ embeds in $L_p$ if and only if the Fourier transform of $\|\cdot\|^p$ is a
positive distribution in the whole $\R^n.$ These characterizations have proved to be useful in the
study of subspaces of $L_p$ and intersection bodies, see \cite[Chapter 6]{K5}. In this section we
prove a similar characterization of spaces that embed in $L_0.$

\begin{Thm}\label{Thm:embed}Let $K$ be an origin symmetric star body in $\mathbb{R}^n$.
The space $(\mathbb{R}^n, \|\cdot\|_K)$ embeds in $L_0$ if and only if the Fourier transform of
$\ln \|x\|_K$ is a negative distribution outside of the origin in $\mathbb{R}^n$.
\end{Thm}
\noindent{\bf Proof.} First, assume that $(\mathbb{R}^n, \|\cdot\|_K)$ embeds in $L_0.$ Let $\phi$
be a non-negative even test function with compact support outside of the origin. By the definition
of embedding in $L_0$, formula (\ref{radon-fourier})(note that $\hat{\hat{\phi}}=(2\pi)^n \phi$
for even $\phi$) and the Fubini theorem,
\begin{eqnarray}\label{eqn:ln^}
\langle \left(\ln \|x\| \right)^\wedge, \phi\rangle &=&\langle \ln \|x\| ,\hat{\phi}(x)
\rangle \nonumber \\
&=& \int_{S^{n-1}}\int_{\mathbb{R}^n} \ln| ( x,\xi)|
\hat{\phi}(x)\ dx\ d\mu(\xi)+
C\int_{\mathbb{R}^n}\hat{\phi}(x) dx \nonumber \\
&=& \int_{S^{n-1}}\langle\ln|t|,\int_{( x,\xi) =t} \hat{\phi}(x)\ dx\rangle\ d\mu(\xi)\nonumber \\
&=& (2\pi)^n \int_{S^{n-1}}\int_{\mathbb{R}}\left(\ln|z|\right)^\wedge(t)\phi(t\xi)\ dt\ d\mu(\xi)
\end{eqnarray}
since $\int_{\mathbb{R}^n}\hat{\phi}(z)dz=(2\pi)^n\phi(0)=0$. Now, the formula for the Fourier
transform of $\ln|t|$ from \cite[p.362]{GS} implies that
\begin{eqnarray}\label{FT:ln|t|}
\left(\ln|z|\right)^\wedge(t)= -\pi |t|^{-1}<0
\end{eqnarray}
outside of the origin, so (\ref{eqn:ln^}) is negative (recall that $\phi$ is non-negative with
support outside of the origin). This means that $\left(\ln \|x\| \right)^\wedge $ is a negative
distribution .

To prove the other direction, note that, by \cite[Section 2.6]{K5}, a distribution that is positive
outside of the origin coincides with a finite Borel measure on every set of the form
$$A\times(a,b)=\left\{ x\in \mathbb{R}^n: x=t\theta, t\in(a,b),\theta \in A\right\}, $$
where $A$ is an open subset of $S^{n-1}$ and $0<a<b<\infty.$

Denote by $\mu=-(\ln\|x\|)^\wedge$. This distribution coincides with a finite Borel measure on each
set $A\times(a,b)$, as above, so for any test function $\phi$ supported outside of the origin

\begin{eqnarray}
\langle -(\ln\|x\|)^\wedge, \phi \rangle &=&\langle \mu, \phi \rangle \nonumber \\
 &=& \int_{\mathbb{R}^n}\phi(x) d\mu(x).
\end{eqnarray}

Now for every test function $\phi$ with support outside of the origin and $t>0$, we have
$\left(\phi(x/t)\right)^\wedge(z)=t^n \hat{\phi}(tz) $, so

\begin{eqnarray}\label{eqn:homog}
\langle \mu(x),\phi(x/t)\rangle &=& - \langle \left(\ln\|x\|\right)^\wedge(x),\phi(x/t)\rangle \nonumber\\
&=& -\int_{\mathbb{R}^n}\ln\|z\|\hat{\phi}(t z)t^n dz \nonumber\\
&=& -\int_{\mathbb{R}^n}\hat{\phi}(\tilde{x})\ln\|\frac{1}{t}\tilde{x}\|d\tilde{x} \nonumber\\
&=& -\int_{\mathbb{R}^n}\hat{\phi}(\tilde{x})\ln\|\tilde{x}\|d\tilde{x}+
\ln|t|\int_{\mathbb{R}^n}\hat{\phi}(\tilde{x})d\tilde{x} \nonumber \\
&=& -\int_{\mathbb{R}^n}\hat{\phi}(\tilde{x})\ln\|\tilde{x}\|d\tilde{x} \nonumber \\
&=&\langle \mu(x), \phi(x)\rangle.
\end{eqnarray}

Let $\chi_{A\times(a,b)}$ be the indicator of the set $A\times(a,b)$. Approximating
$\chi_{A\times(a,b)}$ by test functions and using {\rm (\ref{eqn:homog})}, we get for any $(a,b)
\subset (0,\infty)$ and $A \subset S^{n-1}$
\begin{eqnarray*}
\mu \left(A\times(a,b)\right)&=&\int_{\mathbb{R}^n}\chi_{A\times(a,b)}(x)d\mu(x)\\
&=&\int_{\mathbb{R}^n}\chi_{A\times(1,b/a)}(x/a)d\mu(x)\\
&=&\int_{\mathbb{R}^n}\chi_{A\times(1,b/a)}(x)d\mu(x)\\
&=&\mu(A\times(1,b/a)).
\end{eqnarray*}
Applying this formula $n$ times,
\begin{eqnarray}\label{eqn:measure}
\mu \left(A\times(1,a^n)\right)=n \mu \left(A\times(1,a)\right)
\end{eqnarray}
for $n \in \mathbb{N}$. Moreover, we can extend formula {\rm (\ref{eqn:measure})} to $n \in
\mathbb{R}$. So, for any $a\in(0,\infty)$, $A \subset S^{n-1}$
$$\mu\left( A\times [1,a]\right)= \mu\left( A\times [1,e^{\ln a}]\right)= \ln a \cdot
\mu\left( A\times [1,e]\right)$$

Now for every $(a,b) \subset (0,\infty)$ and $A \subset S^{n-1}$ we have
\begin{eqnarray*}
\mu \left(A\times(a,b)\right)&=&\mu\left(A\times(1,b/a)\right)\\
&=&\ln\left(\frac{b}{a}\right)\mu\left(A\times(1,e)\right)\\
&=&\left(\ln(b)-\ln(a)\right)\mu(A\times(1,e)).
\end{eqnarray*}

Define a measure $\mu_0$ on $S^{n-1}$ by
$$\mu_0(A)=\frac{\mu\left(A\times(a,b)\right)}{\left(\ln(b)-\ln(a)\right)}=\mu(A\times(1,e))$$
for every Borel set $A\subset S^{n-1}$. We have

\begin{eqnarray} \label{eqn:repr}
\int_{S^{n-1}}d\mu_0(\theta)\int_{0}^{\infty}|t|^{-1}\chi_{A\times(a,b)}(t\theta)dt&=
&\left(\ln(b)-\ln(a)\right)\mu_0(A)\nonumber\\
&=&\mu(A\times(a,b)) \nonumber\\
&=&\int_{\mathbb{R}^n}\chi_{A\times(a,b)}(x)d\mu(x)
\end{eqnarray}

Therefore, for an arbitrary even test function $\phi$ supported outside of the origin,

\begin{eqnarray}\label{eqn:mu}
\frac12\langle
\mu,\phi\rangle&=&\int_{S^{n-1}}d\mu_0(\theta)\int_{0}^{\infty}|t|^{-1}\phi(t\theta)dt\nonumber\\
&=&\frac{1}{2} \int_{S^{n-1}}d\mu_0(\theta)\int_{\mathbb{R}}|t|^{-1}\phi(t\theta)dt
\end{eqnarray}
since $A,a,b$ are arbitrary in {\rm (\ref{eqn:repr})}.

Using $\mu=-\left( \ln\|x\|\right)^\wedge$, we get
$$\langle \left( \ln\|x\|\right)^\wedge(\xi),\phi\rangle=
-\int_{S^{n-1}}d\mu_0(\theta)\int_{\mathbb{R}}|t|^{-1}\phi(t\theta)dt.$$

Define a new measure $\tilde{\mu}_0=(2\pi)^{n}\mu_0$. By (\ref{FT:ln|t|}),(\ref{eqn:mu}) and the
connection between the Fourier and Radon transforms
\begin{eqnarray*}
\langle \ln\|x\|,\hat{\phi}(x)\rangle&=&
-\frac{1}{(2\pi)^{n}}\int_{S^{n-1}}d\tilde{\mu}_0(\theta)\int_{\mathbb{R}}|t|^{-1}\phi(t\theta)dt\\
&=&\int_{S^{n-1}}\langle \ln|z|, \mathcal{R}\hat{\phi}(\theta;z)\rangle d\tilde{\mu}_0(\theta)\\
&=&\int_{S^{n-1}}d\tilde{\mu}_0(\theta)\int_{\mathbb{R}}\ln|z|\left( \int_{( x,\theta)=z}\hat{\phi}(x)d x\right) d z\\
&=&\int_{S^{n-1}}d\tilde{\mu}_0(\theta)\int_{\mathbb{R}^n}\ln\left|( x,\theta)
\right|\hat{\phi}(x)d x
\end{eqnarray*}
Thus, we have proved that for any even test function $\phi$ supported outside of the origin
$$\langle (\ln\|x\|)^\wedge,{\phi}\rangle=\Big\langle \left(\int_{S^{n-1}}\ln\left|( x,\theta)
\right|d\tilde{\mu}_0(\theta)\right)^\wedge,{\phi}\Big\rangle.
$$
 Therefore the distributions
$\ln\|x\|$ and $\int_{S^{n-1}}\ln\left|( x,\theta) \right|d\tilde{\mu}_0(\theta)$ can differ only
by a polynomial. Clearly, this polynomial cannot contain terms homogeneous of degree different from
zero, so it is a constant.

 \qed

\begin{Rmk}\label{Rmk:measure mu}
Let $K$ be an infinitely smooth body. From the proof of the previous theorem it follows that the
measure $\mu$ from Definition \ref{Def:L_0} is equal to restriction of the Fourier transform of
$\ln \|x\|_K$ to the sphere. In the next section we are going to prove that this is a function,
therefore

$$d\mu(\xi)=-\frac{1}{(2\pi)^n} \left(\ln \|x\|_K\right)^\wedge(\xi) d\xi.$$
\end{Rmk}

%\begin{Rmk}
In particular, since $\mu$ is a probability measure, for any infinitely smooth body $K$ we get
$$-\frac{1}{(2\pi)^n}\int_{S^{n-1}} \left(\ln \|x\|_K\right)^\wedge(\xi) d\xi=1.$$

%\end{Rmk}

\section{A geometric characterization of subspaces of $L_0$.}

  Let $K$ be an origin symmetric star body
in $\mathbb{R}^n$. The function $\rho_K(x)=\|x\|_K^{-1}$ is called  the {\it radial function} of
$K$. If $x\in S^{n-1}$, $\rho_K(x)$ is the distance from the origin to the boundary of $K$ in the
direction of $x$.

The {\it radial metric} on the set of all origin symmetric star bodies is defined by
$$\rho(K,L)=\max_{x\in S^{n-1}} |\rho_K(x)-\rho_L(x)|.$$

 Let $\xi\in S^{n-1}$ and $(x,\xi) =t$ be the hyperplane orthogonal to $\xi$ at the distance $t$ from the origin. Define
the {\it parallel section function} of a star body $K$ in the direction of $\xi$ by

$$A_{K,\xi}(t)=\mathrm{vol}_{n-1}(K\cap \{( x,\xi) =t\}), \ \ \ t \in \mathbb{R}.$$

Let $f$ be an integrable continuous function on $\mathbb R$, $m$-times continuously differentiable
in some neighborhood of zero, $m \in \mathbb{N}.$ For a number  $q\in (m-1,m)$ the {\it fractional
derivative} of the order $q$ of the function $f$ at zero is defined as follows \cite[Section
2.5]{K5}:

\begin{eqnarray*}
f^{(q)}(0)=\frac{1}{\Gamma(-q)} \int_0^\infty t^{-1-q}\Big( f^{}(t)-f^{}(0)-t f'(0)- \cdots -&& \\
-\frac{t^{m-1}}{(m-1)!}f^{(m-1)}(0)\Big)dt.&&
\end{eqnarray*}

Note, that fractional derivatives of integer orders coincide with usual derivatives up to a sign:
$$f^{(k)}(0)=(-1)^k\frac{d^k}{dt^k}f(t)|_{t=0}.$$

It was shown in \cite{GKS} that if $K$ has an infinitely smooth boundary then the fractional
derivatives of $A_{K,\xi}(t)$ can be computed in terms of the Fourier transform  of the Minkowski
functional raised to certain powers. Namely, for $q\in \mathbb{C}$, $q\ne n-1$,

\begin{eqnarray}\label{eqn:derivA}
A_{K,\xi}^{(q)}(0)=\frac{\cos\frac{q\pi}{2}}{\pi(n-q-1)}\left(\|x\|_K^{-n+q+1}\right)^\wedge (\xi),
\end{eqnarray}
and, in particular, $\left(\|x\|_K^{-n+q+1}\right)^\wedge $ is a continuous function on
$\R^n\setminus \{0\}$. Here we extend $A_{K,\xi}^{(q)}(0)$ from the sphere to the whole $\mathbb{R}^n$ as a
homogeneous function of the variable $\xi$ of degree $-q-1$. Note that $\langle
A_{K,\xi}^{(q)}(0),\phi\rangle$ is an analytic function of $q$ for any fixed test function $\phi$.

In our next Theorem we use a limiting argument to extend formula (\ref{eqn:derivA}) to the case
$q=n-1$.

 Let $\mathcal{D}$ be an open set in $\mathbb{R}^n$, $f,g$ two distributions. We say that
$f=g$ on $\mathcal{D}$ if $\langle f, \phi \rangle=\langle g, \phi \rangle$ for any test function
$\phi$ with compact support in $\mathcal{D}$.

\begin{Thm}\label{Thm:n-1-deriv}
Let $K$ be an infinitely smooth origin symmetric star body in $\mathbb{R}^n$. Extend
$A_{K,\xi}^{(n-1)}(0)$ to a homogeneous function of degree $-n$ of the variable $\xi \in
\mathbb{R}^n \setminus \{0\}$. Then $(\ln \|\cdot\|_K)^\wedge$ is a continuous function
on $\R^n\setminus \{0\}$ and

\begin{eqnarray}\label{eqn:A^{(n-1)}}
A_{K,\xi}^{(n-1)}(0)=-\frac{\cos(\pi (n-1)/2)}{\pi} \left( \ln \|\cdot\|_K\right)^\wedge (\xi),
\end{eqnarray}
as distributions (of the variable $\xi$) acting on test functions with compact support outside of
the origin. In particular,

i) if $n$ is odd

$$\left( \ln \|x\|_K\right)^\wedge (\xi)=(-1)^{(n+1)/2}\pi A_{K,\xi}^{(n-1)}(0),  \ \ \xi \in \mathbb{R}^n\setminus \{0\}$$

ii) if $n$ is even, then for    $\xi \in \mathbb{R}^n\setminus \{0\}$,

$$\left( \ln \|x\|_K\right)^\wedge (\xi)=a_n\int_0^{\infty}\frac{A_\xi(z)-A_\xi(0)-
A''_\xi(0)\frac{z^2}{2}-...-A^{n-2}_\xi(z)\frac{z^{n-2}}{(n-2)!}}{z^n}dz,$$ where
$a_n=2(-1)^{n/2+1}(n-1)!$

\end{Thm}
\noindent{\bf Proof.} Let us start with the case where $n$ is odd.
 Let $\phi$ be a test function supported outside of the origin.

Using formula (\ref{eqn:derivA}) for $q$ close to $n-1$, we have
\begin{eqnarray*}
\langle A_{K,\xi}^{(q)}(0),\phi(\xi)\rangle&=&\frac{\cos(\pi
q/2)}{\pi(n-q-1)}\langle\left(\|x\|^{-n+q+1}\right)^\wedge(\xi),\phi(\xi)\rangle\\
&=&\frac{\cos(\pi q/2)}{\pi(n-q-1)}\langle\|x\|^{-n+q+1},\hat{\phi}(x)\rangle\\
&=&\frac{\cos(\pi q/2)}{\pi(n-q-1)}\int_{\mathbb{R}^n}\|x\|^{-n+q+1}\hat{\phi}(x)d x\\
&=&\frac{\cos(\pi q/2)}{\pi(n-q-1)}\int_{\mathbb{R}^n}\left(\|x\|^{-n+q+1}-1\right)\hat{\phi}(x)d x\\
&&+\frac{\cos(\pi q/2)}{\pi(n-q-1)}\int_{\mathbb{R}^n}\hat{\phi}(x)d x\\
&=&\frac{\cos(\pi q/2)}{\pi}\int_{\mathbb{R}^n}\frac{\|x\|^{-n+q+1}-1}{n-q-1}\hat{\phi}(x)d x,
\end{eqnarray*}
since $\int_{\mathbb{R}^n}\hat{\phi}(x)d x=(2\pi)^n\phi(0)=0.  $ Taking the limit of both sides as
$q \to n-1$, we get
$$\langle A_{K,\xi}^{(n-1)}(0),\phi(\xi)\rangle=\Big\langle -\frac{\cos(\pi (n-1)/2)}{\pi}\left(\ln\|x\|\right)^\wedge(\xi),\phi(\xi)\Big\rangle$$
since
\begin{eqnarray*}
\lim_{q\to n-1}\int_{\mathbb{R}^n}\frac{\|x\|^{-n+q+1}-1}{n-q-1}\hat{\phi}(x)d x&=&-\int_{\mathbb{R}^n}\ln\|x\|\hat{\phi}(x)d x\\
&=&\langle -\left(\ln\|x\|\right)^\wedge(\xi),\phi(\xi)\rangle.
\end{eqnarray*}

When $n$ is odd the formula of i) follows immediately.

When $n$ is even, both sides of (\ref{eqn:A^{(n-1)}}) are equal to zero, and we repeat the
reasoning from Theorem 1 in \cite{GKS}. Divide both sides of (\ref{eqn:derivA}) by $\cos(\frac{\pi
q}{2})$

$$\Big\langle\frac{\left(\|x\|_K^{-n+q+1}\right)^\wedge(\xi)}{(n-q-1)},\phi(\xi)\Big\rangle=\pi\Big\langle\frac{A_{K,\xi}^{(q)}(0)}{\cos\frac{\pi q}{2}}
,\phi(\xi)\Big\rangle $$ and take the limit of both sides when $q \to n-1$.

We have already proved that
$$\lim_{q\to n-1}\Big{\langle} \frac{\left(\|x\|_K^{-n+q+1}\right)^\wedge(\xi)} {(n-q-1)} ,\phi(\xi)\Big{\rangle} =
\langle -\left(\ln\|x\|\right)^\wedge(\xi),\phi(\xi)\rangle$$ for any test function $\phi$
supported outside of the origin.

To compute the limit of $\displaystyle\frac{A_{K,\xi}^{(q)}(0)}{\cos\frac{q\pi}{2}}$ we use the
definition of fractional derivatives in exactly the same way as it was done in \cite[Theorem
1]{GKS}.

$$\lim_{q \to n-1}\Gamma(-q) A_{K,\xi}^{(q)}(0)=\int_0^{\infty}\frac{A_\xi(z)-A_\xi(0)-
A''_\xi(0)\frac{z^2}{2}-...-A^{n-2}_\xi(z)\frac{z^{n-2}}{(n-2)!}}{z^n}dz$$

and

$$\lim_{q \to n-1}\Gamma(-q)\sin\frac{(q+1)\pi}{2}=\frac{\pi}{2}(-1)^{n/2}\frac{1}{(n-1)!}. $$

Combining these two formulas we get the formula in the statement ii) of the Theorem. \qed

\medbreak

An immediate application of Theorem \ref{Thm:n-1-deriv} is

\begin{Cor}
Let $K$ be an infinitely smooth body in $\mathbb{R}^n$. Then

i) if $n$ is odd, $(\mathbb{R}^n, \|\cdot\|_K)$ embeds in $L_0$ if and only if
$$ (-1)^{(n-1)/2} A_{K,\xi}^{(n-1)}(0)\ge 0, \quad \forall \xi\in S^{n-1}; $$%\quad \mbox{ if }n \mbox{ is odd,}$$

ii) if $n$ is even, $(\mathbb{R}^n, \|\cdot\|_K)$ embeds in $L_0$ if and only if, for every $
\xi\in S^{n-1}$,
$$(-1)^{(n+2)/2}\int_0^{\infty}\frac{A_\xi(z)-A_\xi(0)-
A''_\xi(0)\frac{z^2}{2}-...-A^{n-2}_\xi(z)\frac{z^{n-2}}{(n-2)!}}{z^n}dz\ge 0.$$

\end{Cor}

\begin{Cor}\label{Cor:3dim}
Every 3-dimensional normed space $(\mathbb{R}^n, \|\cdot\|_K)$ embeds in $L_0$.
\end{Cor}
\noindent {\bf Proof.}  The unit ball $K$ of a normed space is an origin-symmetric convex body.
First assume that $K$ is infinitely smooth. By Brunn's theorem the central section of a convex body
has maximal volume among all sections perpendicular to a given direction. Therefore, for any $\xi$
the function $A_{K,\xi}(t)$ attains its maximum at $t=0$, hence $A_{K,\xi}''(0)\le 0$. So, by
Theorem \ref{Thm:n-1-deriv}, for smooth convex bodies in $\mathbb{R}^3$ the distribution
$-(\ln\|x\|)^\wedge $ is positive outside of the origin, and our result follows from Theorem
\ref{Thm:embed}. For general convex bodies the result
follows from the facts that any convex body can be approximated by smooth convex bodies and that
positive definiteness is preserved under limits. In fact, let $\{K_i\}$ be a sequence of infinitely
smooth convex bodies that approach $K$ in the radial metric. Then for any non-negative test
function $\phi$ supported outside of the origin we have
$$ - \int_{\mathbb{R}^n}\ln \|x\|_{K_i} \hat{\phi}(x)dx=\langle -\ln\|x\|_{K_i}, \hat{\phi}(x) \rangle = \langle -(\ln\|x\|_{K_i})^\wedge (\xi), \phi
(\xi)\rangle\ge 0$$

Since  $K_i$ approximate $K$ there is a constant $C>0$, such that $$\left|\ln \|x\|_{K_i}\right|\le
C+ \left|\ln |x|_2\right|,$$ therefore the functions $|\ln \|x\|_{K_i} \hat{\phi}(x)|$ are
majorated by an integrable function $(C+|\ln |x|_2 |) |\hat{\phi}(x)|$ and by the Lebesgue
Dominated Convergence Theorem we get
\begin{eqnarray*}-\lim_{i\to \infty}\int_{\mathbb{R}^n}\ln \|x\|_{K_i}
\hat{\phi}(x)dx &=& -\int_{\mathbb{R}^n}\ln \|x\|_{K}
\hat{\phi}(x)dx\\
&=&\langle -(\ln\|x\|_K)^\wedge (\xi), \phi (\xi)\rangle\ge 0
\end{eqnarray*}

\qed

Our next result shows that that the previous statement is no longer true in $\mathbb{R}^n$, $n\ge
4$.

\begin{Thm}
There exists an origin-symmetric convex body $K$ in $\mathbb{R}^n$, $n\ge 4$ so that the space
$(\mathbb{R}^n, \|\cdot\|_K)$ does not embed in $L_0$.
\end{Thm}
\noindent{\bf Proof.} It is enough to construct a convex body for which the distribution $-(\ln
\|x\|)^\wedge$ is not positive. The construction will be similar to that from \cite{GKS}.

Define $f_N(x)=(1-x^2-Nx^4)^{1/3}$, let $a_N>0$ be such that $f_N(a_N)=0$ and $f_N(x)>0$ on the
interval $(0,a_N)$. Define a body $K$ in $\mathbb{R}^4$ by

$$K=\{(x_1,x_2,x_3,x_4)\in \mathbb{R}^4: \, x_4\in [-a_N,a_N]\,  \mathrm{ and }\,
\sqrt{x_1^2+x_2^2+x_3^2}\le f_N(x_4) \}.$$ The body $K$ is strictly convex and infinitely smooth.
By Theorem \ref{Thm:n-1-deriv},
$$-\left( \ln \|x\|_K\right)^\wedge (\xi)=12\int_0^{\infty}\frac{A_\xi(z)-A_\xi(0)-
A''_\xi(0)\frac{z^2}{2}}{z^4}dz.$$

The function $A_{K,\xi}$ can easily be computed: $$A_{K,\xi}(x)=\frac{4\pi}{3} (1-x^2-Nx^4).$$ We
have
$$\int_0^{\infty}\frac{A_\xi(z)-A_\xi(0)-
A''_\xi(0)\frac{z^2}{2}}{z^4}dz=\frac{4\pi}{3}(-Na_N+\frac{1}{a_N}-\frac{1}{3a_N^3}).$$ The latter
is  negative for $N$ large enough, because $N^{1/4}\cdot a_N\to 1$ as $N\to \infty$.

\qed

\section{ Addition in $L_0$ }

It is clear from the definition that the class of bodies $K$ for which $(\mathbb{R}^n,
\|\cdot\|_K)$ embeds in $L_0$ is closed with respect to multiplicative summation, i.e. if two
spaces $(\mathbb{R}^n, \|\cdot\|_{K_1})$ and $(\mathbb{R}^n, \|\cdot\|_{K_2})$ embed in $L_0$ and
$K=K_1+_0 K_2$, then $(\mathbb{R}^n, \|\cdot\|_K)$ embeds in $L_0$. In this section we are going to
prove that the unit ball of every space $(\mathbb{R}^n, \|\cdot\|_K)$ that embeds in $L_0$ can be
obtained from the Euclidean ball by means of multiplicative summation, linear transformations and
closure in the radial metric, i.e. it can be approximated in the radial metric by multiplicative
sums of ellipsoids.

Consider the set of bodies $K$ for which $(\mathbb{R}^n, \|\cdot\|_{K})$ embeds in $L_0$.
As mentioned above,
this set is closed with respect to multiplicative summation, also from the proof of Corollary
\ref{Cor:3dim} it follows that this set is closed with respect to limits in the radial metric. Let
us show that it is closed with respect to linear transformations. Suppose that $(\mathbb{R}^n,
\|\cdot\|_{K})$ embeds in $L_0$. By Theorem \ref{Thm:embed} $\left(\ln\|x\|_K\right)^\wedge$ is a
negative distribution outside of the origin. Let $T$ be a linear transformation in $\mathbb{R}^n$,
then for any non-negative test function $\phi$ with support outside of the origin, we have

\begin{eqnarray*}
\langle \left(\ln\|Tx\|_K\right)^\wedge, \phi\rangle&=&\langle \ln\|Tx\|_K, \hat{\phi}(x)\rangle\\
&=&\int_{\mathbb{R}^n} \ln\|Tx\|_K \hat{\phi}(x) dx\\
&=&|\det T|^{-1} \int_{\mathbb{R}^n} \ln\|x\|_K \hat{\phi}(T^{-1}x) dx\\
&=& \int_{\mathbb{R}^n} \ln\|x\|_K \left({\phi}(T^{*}y)\right)^\wedge(x) dx\\
&=& \langle \ln\|x\|_K, \left({\phi}(T^{*}y)\right)^\wedge(x) \rangle,\\
&=& \langle \left(\ln\|x\|_K\right)^\wedge(y), \phi(T^{*}y) \rangle\le 0.
\end{eqnarray*}
So $\left(\ln\|Tx\|_K\right)^\wedge$ is a negative distribution outside of the origin. By Theorem
\ref{Thm:embed}, $(\mathbb{R}^n, \|\cdot\|_{TK})$  embeds in $L_0$.

Moreover, if $\left(\ln\|x\|\right)^\wedge$ is a function, then
\begin{equation}\label{linear transform}
\left(\ln\|Tx\|\right)^\wedge(y)=|\det T|^{-1}\left(\ln\|x\|\right)^\wedge((T^{*})^{-1} y).
\end{equation}

To prove the main result of this section we need a few lemmas.
%{\bf Step 1.}
For a fixed $x\in S^{n-1}$, let $E_{a,b}(x)$ be an ellipsoid with the norm
$$\|\theta\|_{E_{a,b}(x)}=\left(\frac{( x, \theta )^2}{a^2}+ \frac{1-( x, \theta )^2}{b^2}\right)^{1/2}, \quad \mbox{ for } \theta \in S^{n-1}.$$

\begin{Lem}\label{Lem:FTellipsoid} For all $\theta\in S^{n-1}$,
$$\left(\ln  \|\xi\|_{E_{a,b}(x)}\right)_\xi^\wedge(\theta) = -\frac{2^{n-1} \pi^{n/2}\Gamma(n/2)}{a^{n-1}b} \|\theta\|_{E_{b,a}(x)}^{-n}. $$
\end{Lem}
\noindent{\bf Proof.} For $-n<\lambda<0$ the following formula holds (see \cite[p.192]{GS}):

$$\left(|x|_2^\lambda\right)^\wedge(\xi)=2^{\lambda+n}\pi^{n/2}\frac{\Gamma((\lambda+n)/2)}{\Gamma(-\lambda/2)}|\xi|_2^{-\lambda-n}.$$
Dividing both sides by $\lambda$, using the formula $x\Gamma(x)=\Gamma(1+x)$  and sending
$\lambda\to 0$ we get

\begin{eqnarray*}
\left(\ln |x|_2\right)^\wedge(\xi)=-2^{n-1}\pi^{n/2} \Gamma(n/2)|\xi|_2^{-n},
\end{eqnarray*}
as distributions outside of the origin. Note that, by rotation, it is enough to prove  Lemma for
the ellipsoids $E_{a,b}(x)$ with $x=(0,0,\dots,0,1)$.
$$\|\xi\|_{E_{a,b}(x)}=\left(\frac{\xi_n^2}{a^2}+ \frac{\xi_1^2+\cdots+\xi_{n-1}^2}{b^2}\right)^{1/2}.$$ %\quad \mbox{ where } \xi,x\in S^{n-1}$$
Since this norm can be obtained from the Euclidean norm by an obvious linear transformation, one
can use formula (\ref{linear transform}) to get
\begin{eqnarray*}
\left(\ln  \|\xi\|_{E_{a,b}(x)}\right)_\xi^\wedge(\theta) &=& -2^{n-1}
\pi^{n/2}\Gamma(n/2){ab^{n-1}} \|\theta\|_{E_{1/a,1/b}(x)}^{-n}\\
&=&-\frac{2^{n-1} \pi^{n/2}\Gamma(n/2)}{a^{n-1}b} \|\theta\|_{E_{b,a}(x)}^{-n}.
\end{eqnarray*}
 \qed

\begin{Lem}\label{Step1}
Let $K$ be a star body, then $\ln\|x\|_K$ can be  approximated in the space $C(S^{n-1})$ by the
functions of the form
\begin{equation}\label{eqn:delta-seq}
f_{a,b}(x)=\frac{1}{|S^{n-1}|a^{n-1}b} \int_{S^{n-1}}\ln\|\theta\|_K \|\theta\|_{E_{b,a}(x)}^{-n} d\theta, %\int_{S^{n-1}}\ln  \|\xi \|_{E_{a,b}(x)} d\mu(\xi)
\end{equation} as $a\to 0$ and $b$ is fixed.
\end{Lem}

\noindent{\bf Proof.} The proof is similar to that of \cite[Lemma 2]{GW}. First, note that the
space $\mathbb{R}^n$ with the Euclidean norm embeds in $L_0$, so $(\mathbb{R}^n, \|\cdot\|_E)$
embeds in $L_0$ for any ellipsoid $E$ with center at the origin. Therefore, by Remark
\ref{Rmk:measure mu} and Lemma \ref{Lem:FTellipsoid} we get

$$\int_{S^{n-1}}\frac{ 1}{|S^{n-1}| a^{n-1}b} \|\theta\|_{E_{b,a}(x)}^{-n} d\theta =1,$$
for all values of $a$ and $b$. From now on $b$ will be fixed.

We have

\begin{eqnarray*}
&&\left|\ln\|x\|_K- \frac{1}{|S^{n-1}|a^{n-1}b} \int_{S^{n-1}}\ln\|\theta\|_K
\|\theta\|_{E_{b,a}(x)}^{-n} d\theta\right|\\
&&\le \frac{1}{|S^{n-1}| a^{n-1}b} \int_{S^{n-1}}\Big|\ln\|x\|_K - \ln\|\theta\|_K
\Big| \|\theta\|_{E_{b,a}(x)}^{-n} d\theta\\
&&= \frac{1}{|S^{n-1}| a^{n-1}b} \int_{|(x,\theta)|\ge\delta}\Big|\ln\|x\|_K - \ln\|\theta\|_K
\Big| \|\theta\|_{E_{b,a}(x)}^{-n} d\theta\\
&&+ \frac{1}{|S^{n-1}|a^{n-1}b} \int_{|(x,\theta)|<\delta}\Big|\ln\|x\|_K - \ln\|\theta\|_K
\Big| \|\theta\|_{E_{b,a}(x)}^{-n} d\theta\\
&&=I_1+I_2.
\end{eqnarray*}
For the first integral $I_1$ use the uniform continuity of $\ln\|x\|_K$ on the sphere. For any
given $\epsilon>0$ there exists $\delta\in(0,1)$, $\delta$ close to $1$, so that
$|(x,\theta)|\ge\delta$ implies $\Big|\ln\|x\|_K - \ln\|\theta\|_K \Big|<\epsilon/2$. Therefore
\begin{eqnarray*}
I_1&=& \frac{1}{|S^{n-1}| a^{n-1}b} \int_{|(x,\theta)|\ge\delta}\Big|\ln\|x\|_K - \ln\|\theta\|_K
\Big| \|\theta\|_{E_{a,b}(x)}^{-n} d\theta\\
&\le&\frac{\epsilon}{2} \left[\frac{1}{|S^{n-1}|a^{n-1}b} \int_{|(x,\theta)|\ge\delta}
\|\theta\|_{E_{a,b}(x)}^{-n} d\theta \right]\le \frac{\epsilon}{2}.
\end{eqnarray*}

Now fix $\delta$ chosen above and estimate the integral $I_2$ as follows
\begin{eqnarray*}
I_2&=& \frac{1}{ |S^{n-1}| a^{n-1}b} \int_{|(x,\theta)|<\delta}\Big|\ln\|x\|_K - \ln\|\theta\|_K
\Big| \|\theta\|_{E_{b,a}(x)}^{-n} d\theta\\
&\le&  \frac{ C(n,b,K)}{a^{n-1}}
\int_{|(x,\theta)|<\delta} \|\theta\|_{E_{b,a}(x)}^{-n} d\theta,\\
\end{eqnarray*}
where $$C(n,b,K)=\frac{2 \max_{S^{n-1}}|\ln\|x\|_K|}{  |S^{n-1}|^{} b }.$$ For the latter integral
we use an elementary formula (see e.g. \cite[Section 6.4]{K5})

$$\int_{|(x,\theta)|<\delta}f((x,\theta))d\theta= |S^{n-2}| \int_{-\delta}^{\delta} (1-t^2)^{(n-3)/2}
f(t)dt, \quad \mbox{ for } x\in S^{n-1}.$$ Now,
\begin{eqnarray*}
I_2&\le&   \frac{ C(n,b,K) |S^{n-2}| }{a^{n-1}} \int_{-\delta}^{\delta} (1-t^2)^{(n-3)/2}
\left(\frac{t^2}{b^2}+
    \frac{1-t^2}{a^2}\right)^{-n/2}dt\\
&\le&   \frac{ C(n,b,K) |S^{n-2}| }{a^{n-1}} \int_{-\delta}^{\delta} (1-t^2)^{(n-3)/2} \left(
    \frac{1-t^2}{a^2}\right)^{-n/2}dt\\
    &=&a\cdot  C(n,b,K) |S^{n-2}|  \int_{-\delta}^{\delta} (1-t^2)^{-3/2}dt\\
    &\le& a\cdot C(n,b,K) |S^{n-2}| \frac{2 \delta} {(1-\delta^2)^{3/2}}.
\end{eqnarray*}
Now we can choose $a$ so small that $I_2\le \epsilon/2$. \qed

\begin{Lem}\label{Step2}
If $\mu$ is a probability measure on $S^{n-1}$ and $a,b>0$, then the function

$$f(x)=\int_{S^{n-1}} \ln  \|\xi \|_{E_{a,b}(x)} d\mu(\xi)$$ can be  approximated in $C(S^{n-1})$ by the sums of the form
$$\sum_{i=1}^m \frac{1}{p_i} \ln \|x\|_{E_i},$$
where $E_1$,...,$E_m$ are ellipsoids and $1/p_1+\cdots +1/p_m=1$.
\end{Lem}
\noindent{\bf Proof.} Let $\sigma>0$ be a small number and choose a finite covering of the sphere
by spherical $\sigma$-balls $B_\sigma (\eta_i)=\{\eta\in S^{n-1}: |\eta-\eta_i|<\sigma\}$,
$\eta_i\in S^{n-1}$, $i=1,\dots, m=m(\delta)$. Define
$$\widetilde{B}_\sigma(\xi_1)=B_\sigma (\xi_1)$$
and
$$\widetilde{B}_\sigma(\xi_i)=B_\sigma (\xi_i)\setminus\bigcup_{j=1}^{i-1}B_\sigma(\xi_j), \quad\mbox{ for } i=2,...,m.$$

Let $1/p_i=\mu(\widetilde{B}_\sigma(\xi_i))$. Clearly, $1/p_1+\cdots+1/p_m=1$.

Let $\rho(E_{a,b}(\xi),x)$ be the  radial function of the ellipsoid $E_{a,b}(\xi)$, that is
$$\rho(E_{a,b}(\xi),x)=\|x\|_{E_{a,b}(\xi)}^{-1}.$$

Note that $\rho(E_{a,b}(\xi),x)=\rho(E_{a,b}(x),\xi),$ therefore
$$|\rho(E_{a,b}(\xi),x)- \rho(E_{a,b}(\theta),x)|\le C_{a,b} |\xi - \theta|,$$
with a constant $C_{a,b}$ that depends  on $a$ and $b$. Also note that, since we consider $a$ close
to zero and $b$ fixed, we may assume
$$a\le \rho(E_{a,b}(\xi),x)\le b,\qquad x\in S^{n-1}.$$
Then,
\begin{eqnarray*}
\left|\int_{S^{n-1}} \ln \rho(E_{a,b}(\xi),x) d\mu(\xi)-\sum_{i=1}^m\frac{1}{p_i}\ln \rho(E_{a,b}(\xi_i),x)\right|=\\
=\left|\sum_{i=1}^m\left(\int_{\widetilde{B}_\sigma(\xi_i)} \ln \rho(E_{a,b}(\xi),x)
d\mu(\xi)-\int_{\widetilde{B}_\sigma(\xi_i)} \ln \rho(E_{a,b}(\xi_i),x) d\mu(\xi)\right)\right|\le\\
\le \sum_{i=1}^m \int_{\widetilde{B}_\sigma(\xi_i)} \left| \ln\frac{
\rho(E_{a,b}(\xi),x)}{\rho(E_{a,b}(\xi_i),x)} \right|d\mu(\xi)\le\\
\le \sum_{i=1}^m \int_{\widetilde{B}_\sigma(\xi_i)} \left| \ln\frac{\rho(E_{a,b}(\xi_i),x)+
[\rho(E_{a,b}(\xi),x)-\rho(E_{a,b}(\xi_i),x)]}{\rho(E_{a,b}(\xi_i),x)} \right|d\mu(\xi)\le\\
\le \sum_{i=1}^m \int_{\widetilde{B}_\sigma(\xi_i)} \left|\ln (1\pm C'_{a,b} |\xi-\xi_i|)\right|
d\mu(\xi)\le\\
\le  \left|\ln (1\pm C'_{a,b} \sigma)\right|,
\end{eqnarray*}
and the result follows since $\sigma$ is arbitrarily small.

 \qed

Now we are ready to prove the following

\begin{Thm}\label{Thm:approximation}
Let $K$ be an origin symmetric star body in $\mathbb{R}^n$. The space $(\mathbb{R}^n, \|\cdot\|_K)$
embeds in $L_0$ if and only if $ \|x\|_K$ is the limit (in the radial metric) of finite products
$\|x\|_{E_1}^{1/p_1}\cdots \|x\|_{E_m}^{1/p_m}$, where $E_1$,...,$E_m$ are ellipsoids and
$1/p_1+\cdots +1/p_m=1$.
\end{Thm}
\noindent{\bf Proof.} The ``if" part is a consequence of the fact that $L_0$ is closed with respect
to the three operations as discussed above.

The proof of ``only if" part  easily  follows the Lemmas we have proved.

Suppose that $(\mathbb{R}^n, \|\cdot\|_{K})$ embeds in $L_0$ with the corresponding probability
measure $\mu$ on $S^{n-1}$ and constant $C$. By Remark \ref{Rmk:measure mu}, $(\mathbb{R}^n,
\|\cdot\|_{E_{a,b}(x)})$ embeds in $L_0$ with the measure $-\frac{1}{(2\pi)^n} \left(\ln
\|x\|_E\right)^\wedge(\theta) d\theta$ and some constant $C_{E_{a,b}}$. Note, this constant does
not depend on $x$. We have

\begin{eqnarray*}
&&\int_{S^{n-1}}\ln  \|\xi \|_{E_{a,b}(x)} d\mu(\xi)\\
&=& \int_{S^{n-1}} \int_{S^{n-1}} \ln |(\xi, \theta ) |
\left(-\frac{1}{(2\pi)^n}\right) \left(\ln \|x\|_{E_{a,b}(x)}\right)^\wedge(\theta) d\theta d\mu(\xi) +C_{E_{a,b}}\\
&=& \int_{S^{n-1}} \left[\int_{S^{n-1}} \ln |(\xi, \theta )| d\mu(\xi)+C_K\right]
\left(-\frac{1}{(2\pi)^n}\right) \left(\ln \|x\|_{E_{a,b}(x)}\right)^\wedge(\theta) d\theta  \\
&& + C_{E_{a,b}}-C_K\\
&=& \int_{S^{n-1}} \ln\|\theta\|_K
\left(-\frac{1}{(2\pi)^n}\right) \left(\ln \|x\|_{E_{a,b}(x)}\right)^\wedge(\theta) d\theta  +C_{E_{a,b}}-C_K\\
&=& \int_{S^{n-1}} \ln\|\theta\|_K
\left(-\frac{1}{(2\pi)^n}\right) \left(\ln \|x\|_{E_{a,b}(x)}\right)^\wedge(\theta) d\theta  +C_{E_{a,b}}-C_K\\
&=& \frac{1}{|S^{n-1}| a^{n-1}b} \int_{S^{n-1}}\ln\|\theta\|_K \|\theta\|_{E_{b,a}(x)}^{-n} d\theta
+C_{E_{a,b}}-C_K
\end{eqnarray*}

In Lemma \ref{Step1} we proved that $\ln\|x\|_K$ can be uniformly approximated by the integrals of
the form
\begin{equation*}
\frac{1}{|S^{n-1}|a^{n-1}b} \int_{S^{n-1}}\ln\|\theta\|_K \|\theta\|_{E_{b,a}(x)}^{-n} d\theta, %\int_{S^{n-1}}\ln  \|\xi \|_{E_{a,b}(x)} d\mu(\xi)
\end{equation*} as $a\to 0$.
Therefore, using the previous calculations, one can see that $\ln\|x\|_K$ can be uniformly
approximated by $$\int_{S^{n-1}}\ln  \|\xi \|_{E_{a,b}(x)} d\mu(\xi)+C'.$$

%{\bf Step 2.}

 Hence, by Lemma \ref{Step2},  $\ln\|x\|_K$ can be uniformly
approximated by the sums $$\sum_{i=1}^m \frac{1}{p_i} \ln \|x\|_{E_i}+C'.$$

Replacing $E_1$ by another ellipsoid $E_1'$ given by $\|x\|_{E'_1}^{1/p_1}=e^{C'}
\|x\|_{E_1}^{1/p_1}$, we get the statement of the Theorem.

\qed

\begin{Cor}
Any convex body in $ \mathbb{R}^3$ can be obtained from the Euclidean unit ball by means of
three operations: linear transformations, multiplicative addition and closure in the radial metric.
\end{Cor}

\noindent{\bf Proof.} %The result follows from Corollary \ref{Cor:3dim}, Theorem \ref{Thm:approximation} and the fact that $L_0$ is closed...
 As was proved in Theorem \ref{Thm:approximation}, any convex body can be approximated by
the finite products of the type $\|x\|_{E_1}^{1/p_1}\cdots \|x\|_{E_m}^{1/p_m}$. Since any number
$1/p$ can be approximated by the sums $$\frac{1}{2^{i_1}}+\frac{1}{2^{i_2}}+\cdots
+\frac{1}{2^{i_k}},$$ the result follows.

\qed

A  proof similar to that of Theorem \ref{Thm:approximation} can be used to show that the previous
theorem holds for $p$-summation with $-1<p<1$, $p\ne0$,  in place of the multiplicative summation.

\begin{Thm}\label{Thm:approximation^p}
Let $K$ be an origin symmetric star body in $\mathbb{R}^n$. The space $(\mathbb{R}^n, \|\cdot\|_K)$
embeds in $L_p$, $-1<p<1$, $p\ne 0$   if and only if $ \|x\|_K^p $ is the limit (in the radial
topology) of finite sums $\|x\|_{E_1}^{p}+\cdots +\|x\|_{E_m}^{p}$, where $E_1$,...,$E_m$ are
ellipsoids.
\end{Thm}

\section{Confirming the place of $L_0$ in the scale of $L_p$-spaces.}

In this section we establish the relations between embedding in $L_0$ and in $L_p$ with $p\ne0$,
which confirm the place of $L_0$ between $L_p$ with $p>0$ and $p<0$. We are going to use the
following result from \cite[Theorem 1]{K3}:

\begin{Thm} \label{Thm:embed-p}
An n-dimensional homogeneous space $(\mathbb{R}^n, \|\cdot\|_K)$ embeds in $L_{-p}$,
$p\in (0,n)$ if and only if $\|x\|_K^{-p}$ is a positive definite distribution.
\end{Thm}

We also use a well-known result of P.Levy (see \cite[p.189]{BL}, also  \cite{BDK} for the infinite
dimensional case):

\begin{Thm}\label{Thm:Levy}
 A space $(\mathbb{R}^n, \|\cdot\|_K)$ embeds in $L_p$, $p\in (0,2]$ if and only if the
function $\exp(-\|x\|_K^p)$ is positive definite.
\end{Thm}

Now we are ready to prove

\begin{Thm}
Let $K$ be an origin symmetric star body in $\mathbb{R}^n$. If the space $(\mathbb{R}^n,
\|\cdot\|_K)$ embeds in $L_0$ then it also embeds in $L_{-p}$, $0<p<n$.
\end{Thm}
\noindent{\bf Proof.} By Theorem \ref{Thm:approximation}, $\|x\|_K$ is the limit of  finite
products $\|x\|_{E_1}^{1/p_1}\cdots \|x\|_{E_m}^{1/p_m}$. Consider $\|x\|_K^{-p}$ for $0<p<n$. It
is the limit of the products of the form $\|x\|_{E_1}^{-p/p_1}\cdots \|x\|_{E_m}^{-p/p_m}$. Using
the formula
$$\|x\|^{-p}=\frac{2}{\Gamma(p/2)}\int_0^\infty t^{p-1} \exp(-t^2\|x\|^2) dt,$$
we get

\begin{eqnarray*}
\|x\|_{E_1}^{-p/p_1}\cdots \|x\|_{E_m}^{-p/p_m}&=&C\int_0^\infty \cdots\int_0^\infty
t_1^{p/p_1-1}\cdots
t_m^{p/p_m-1}\times\\
&&\times\exp(-t_1^2\|x\|_{E_1}^2-\cdots-t_m^2\|x\|_{E_m}^2) dt_1\cdots dt_m,
\end{eqnarray*}
where $$C=\frac{2^m}{\Gamma(p/2p_1)\cdots \Gamma(p/2p_m)}.$$ Therefore, for any non-negative test
function $\phi$ we have
\begin{eqnarray*}
&&\langle (\|x\|_{E_1}^{-p/p_1}\cdots \|x\|_{E_m}^{-p/p_m})^\wedge(\xi),\phi(\xi)\rangle=\langle
\|x\|_{E_1}^{-p/p_1}\cdots \|x\|_{E_m}^{-p/p_m},\hat{\phi}(x)\rangle=\\
&&=C\int_0^\infty \cdots\int_0^\infty t_1^{p/p_1-1}\cdots t_m^{p/p_m-1}\times\\
&&\times \langle \exp(-t_1^2\|x\|_{E_1}^2-\cdots-t_m^2\|x\|_{E_m}^2),
\hat{\phi}(x)\rangle dt_1\cdots dt_m=\\
&&=C\int_0^\infty \cdots\int_0^\infty t_1^{p/p_1-1}\cdots t_m^{p/p_m-1}\times\\
&&\times\langle (\exp(-t_1^2\|x\|_{E_1}^2-\cdots-t_m^2\|x\|_{E_m}^2))^\wedge(\xi), \phi(\xi)\rangle
dt_1\cdots dt_m.
\end{eqnarray*}
We claim that the latter expression is non-negative. Indeed,  $(\mathbb{R}^n,\|x\|_E$) embeds in
$L_2$  for any ellipsoid, therefore the $2$-sum of ellipsoids
$t_1^2\|x\|_{E_1}^2+\cdots+t_m^2\|x\|_{E_m}^2$ embeds in $L_2$, and hence by Theorem
\ref{Thm:Levy}, the function $\exp(-t_1^2\|x\|_{E_1}^2-\cdots-t_m^2\|x\|_{E_m}^2)$ is positive
definite. Now the fact that $\langle ( \|x\|^{-p}_K )^\wedge , \phi\rangle\ge 0$ follows by an
approximation argument, as in Corollary \ref{Cor:3dim}.

 \qed

\begin{Thm}{\label{Thm:less0}}
Let $K$ be an origin symmetric star body in $\mathbb{R}^n$. If the space $(\mathbb{R}^n,
\|\cdot\|_K)$ embeds in $L_{-p}$ for every  $p \in (0,\epsilon)$, then it also embeds in $L_0$.

\end{Thm}
\noindent{\bf Proof.} The space $(\mathbb{R}^n, \|\cdot\|_K)$ embeds in $L_{-p}$, so by Theorem
\ref{Thm:embed-p} the distribution $\|x\|^{-p}$ is positive definite. Then for every non-negative
test function $\phi$ supported outside of the origin,
\begin{eqnarray*}
-\int_{\mathbb{R}^n} \ln\|x\| \hat{\phi}(x) dx&=&\lim_{p\to
0}\frac{1}{p}\int_{\mathbb{R}^n}(\|x\|^{-p}-1) \hat{\phi}(x) dx\\
&=&\lim_{p\to 0}\frac{1}{p}\int_{\mathbb{R}^n}\|x\|^{-p} \hat{\phi}(x) dx\ge0.
\end{eqnarray*}
The result follows from Theorem \ref{Thm:embed}.

\qed

\begin{Thm}
There are normed spaces that embed in $L_0$, but do not embed in $L_p$ for $p>0$.
\end{Thm}

\noindent{\bf Proof.} As  proved above, every $3$-dimensional normed space embeds in $L_0$, hence
$l_q^3$ with $q>2$ does. On the other hand, $l_q^3$, $q>2$ does not embed in $L_p$ for $0<p\le 2$
(see \cite{K1}).

\qed

Let us also mention that one can use the approach of \cite{KK0} to
produce examples in the same spirit.  It follows from \cite{KK0},
Proposition 3.5 that $\mathbb R\oplus_2\ell_1$ does not embed
isometrically into $L_p$ for $p>0$; hence neither does $\mathbb
R\oplus_2\ell_1^n$ for large enough $n.$

\begin{Prop}\label{alt}  For any $n\in\mathbb N$ the space
$\mathbb R\oplus_2\ell_1^n$ embeds in $L_0$.\end{Prop}

 \noindent{\bf Proof.} Let
$(f_n)_{n=1}^{\infty}$ be a sequence of functions on some
probability space which are independent and 1-stable symmetric, so
that $\mathbb E(e^{itf_j})=e^{-|t|}$ (i.e. the $f_j$ have the
Cauchy distribution).  Then it is clear that
$$ \mathbb E \ln|\sum_{j=1}^na_jf_j|= \ln \sum_{j=1}^n|a_j|.$$
Indeed this follows from the fact that
$$ \frac1\pi\int_{-\infty}^{\infty}\frac{\ln|x|}{1+x^2}dx=0.$$
On the other hand if $f=\sum_{j=1}^na_jf_j$ where
$\sum_{j=1}^n|a_j|=1$ then $f$ has the Cauchy distribution and so
has the same distribution as $g_1/g_2$ where $g_1,g_2$ are
independent normalized Gaussians.  Hence \begin{align*} \mathbb
E\ln|a+bf|&= \mathbb E (\ln |ag_2+bg_1|-\ln|g_2|)\\
&= \ln(a^2+b^2)^{\frac12}.\end{align*} Now for any
$a_0,a_1,\ldots,a_n\in\mathbb R$ we have
$$ \mathbb E|a_0+\sum_{j=1}^na_jf_j|= \ln
\left(|a_0|^2+(\sum_{j=1}^n|a_j|)^2\right)^{\frac12}.$$

This shows (using the remarks at the end of \S \ref{sec2}) that
$\mathbb R\oplus_2\ell_1^n$ embeds into $L_0$ for every $n.$\qed

\begin{Thm}
Let $K$ be an origin symmetric star body in $\mathbb{R}^n$. If the space $(\mathbb{R}^n,
\|\cdot\|_K)$ embeds in $L_{p_0}$ , $0<p_0\le 2$, then it also embeds in $L_0$.

\end{Thm}
\noindent{\bf Proof.} Since $(\mathbb{R}^n, \|\cdot\|_K)$ embeds in $L_{p_0}$, $0<p_0\le2$, by
\cite[Theorem 2]{K3} it also embeds in $L_{-p}$ for any $p\in (0,n)$ and hence, by Theorem
\ref{Thm:less0}, it embeds in $L_0$.

 \qed

\end{document}